\documentclass[reqno,10pt,draft]{amsart}

\usepackage[a4paper,left=35mm,right=35mm,top=30mm,bottom=30mm,marginpar=25mm]{geometry} 
\usepackage{amsmath}
\usepackage{esint}
\usepackage{amssymb}
\usepackage{amsthm}
\usepackage{amscd}
\usepackage[ansinew]{inputenc}
\usepackage{color}
\usepackage[final]{graphicx}

\usepackage{rotating}
\usepackage{tikz}
\usetikzlibrary{shapes.geometric, arrows,positioning}
\tikzstyle{model} = [rectangle,  minimum height=2cm,font=\normalsize,text centered, draw=black, fill=black!5]
\tikzstyle{title}=[fill=black!20, text=black,rounded corners,font=\small\bfseries\sffamily]
\usepackage{tikz-cd}

\usepackage{xcolor}
\usepackage{todonotes}
\usepackage{changes}
\definechangesauthor[name={Ana}, color=magenta]{amr}
\definechangesauthor[name={Elvira}, color=blue]{ez}
\definechangesauthor[name={AnaNotes}, color=gray]{toerase}

\usepackage{lscape}

\renewcommand{\epsilon}{\varepsilon}

\numberwithin{equation}{section}

\newtheoremstyle{thmlemcorr}{10pt}{10pt}{\itshape}{}{\bfseries}{.}{10pt}{{\thmname{#1}\thmnumber{ #2}\thmnote{ (#3)}}}
\newtheoremstyle{thmlemcorr*}{10pt}{10pt}{\itshape}{}{\bfseries}{.}\newline{{\thmname{#1}\thmnumber{ #2}\thmnote{ (#3)}}}
\newtheoremstyle{defi}{10pt}{10pt}{\itshape}{}{\bfseries}{.}{10pt}{{\thmname{#1}\thmnumber{ #2}\thmnote{ (#3)}}}
\newtheoremstyle{remexample}{10pt}{10pt}{}{}{\bfseries}{.}{10pt}{{\thmname{#1}\thmnumber{ #2}\thmnote{ (#3)}}}
\newtheoremstyle{ass}{10pt}{10pt}{}{}{\bfseries}{.}{10pt}{{\thmname{#1}\thmnumber{ A#2}\thmnote{ (#3)}}}

\theoremstyle{thmlemcorr}
\newtheorem{theorem}{Theorem}
\numberwithin{theorem}{section}
\newtheorem{lemma}[theorem]{Lemma}

\theoremstyle{thmlemcorr*}
\newtheorem{theorem*}{Theorem}
\newtheorem{lemma*}[theorem]{Lemma}
\newtheorem{corollary*}[theorem]{Corollary}
\newtheorem{proposition*}[theorem]{Proposition}
\newtheorem{problem*}[theorem]{Problem}
\newtheorem{conjecture*}[theorem]{Conjecture}

\theoremstyle{defi}
\newtheorem{definition}[theorem]{Definition}

\theoremstyle{remexample}
\newtheorem{remark}[theorem]{Remark}

\theoremstyle{ass}

\DeclareMathOperator{\esssup}{ess\,sup}

\newcommand{\R}{\mathbb{R}}

\def\XXint#1#2#3{{\setbox0=\hbox{$#1{#2#3}{\int}$} 
		\vcenter{\hbox{$#2#3$}}\kern-.5\wd0}}

\title[Lower semicontinuity of nonlocal supremals]{A sufficient condition for the lower semicontinuity of nonlocal supremal functionals in the vectorial case}
\author{Giuliano Gargiulo}
\address{ Dipartimento di Scienze e Tecnologie,
	Universit\'a degli studi del Sannio,
	via de Sanctis,  82100, Benevento, Italy}
\email{ggargiul@unisannio.it}

\author{Elvira Zappale*}
\address{Dipartimento di Scienze di Base ed Applicate per l'Ingegneria, Sapienza, Universit\'a di Roma, via Antonio Scarpa 16, 00161 Roma, Italy,  *corresponding author}
\email{elvira.zappale@uniroma1.it}
\begin{document}
	\maketitle
\thispagestyle{empty}	
	 \begin{abstract}
In this note we present a sufficient condition ensuring lower semicontinuity for nonlocal supremal functionals of the type $$W^{1,\infty}(\Omega;\mathbb R^d)\ni u \mapsto \esssup_{(x,y)\in \Omega} W(x,y, \nabla u(x),\nabla u(y)),$$ where $\Omega$ is a bounded open subset of $\mathbb R^N$ and $W:\Omega \times \Omega \times \mathbb R^{d \times N}\times \mathbb R^{d \times N} \to \mathbb R$.

		
		
		\noindent\textsc{Date:} \today.
	\end{abstract}




	
	
	
	
	\section{Introduction}
In recent years a great attention has been devoted to nonlocal functionals both in the integral and supremal setting, due to the many applications to peridynamics, machine learning, image processing, etc. \cite{ADK, BCMC1, BCMC2, BCMC3, BMCP, DFKS, HK} and to $L^\infty$ variational problems e.g. \cite{APESAIM, BK, CKM, EP, K, RZ}, among a wide literature.

Motivated by the Direct Methods in the Calculus of Variations the study of necessary and sufficient conditions ensuring lower semicontinuity of such functionals has been conducted in many papers, see \cite{BMC,  Bevan Pedregal,  KRZ, KZ, Pedregal, Pedregal 2, Pedregal3}.	

In particular, given a bounded open set $\Omega \subset \mathbb R^N$, in \cite{KZ} characterizing conditions  for the sequential lower semicontinuity  in $L^\infty(\Omega;\mathbb R^d)$ of the functional $G: v \in L^\infty(\Omega ;\mathbb R^d) \to \mathbb R$ defined as \begin{align}\label{Fesssup}
	  G(v):=\esssup_{(x,y)\in \Omega\times \Omega} W(v(x), v(y))
\end{align} 
have been provided. Furthermore in \cite{KRZ}  necessary and sufficient assumptions on the supremand $W$ have been determined to ensure that, in absence of lower semicontinuity, the sequentially weakly * lower semicontinuous envelope of  $G$ has the same form, i.e. it can be expressed as a double supremal functional. We also emphasize that \cite{KRZ} contains a power-law approximation result for functionals as in \eqref{Fesssup}, which, on the other hand, also appear in their inhomogeneous version in the context of image denoising (cf. \cite{DFKS}).

Unfortunately, analogous results are not available in the context where the fields $v$ satisfy some differential constraint, in particular when $v(x)= \nabla u (x)$, with $u \in W^{1,\infty}(\Omega;\mathbb R^d)$.
In this paper we will show that a sufficient condition for the functional 
$$W^{1,\infty}(\Omega;\mathbb R^d)\ni u \mapsto \esssup_{(x,y)\in \Omega \times \Omega} W(x,y, \nabla u(x),\nabla u(y)),$$
to be weakly * sequentially lower semicontinuous in $W^{1,\infty}(\Omega;\mathbb R^{d \times N})$ is the {\it separate curl Young quasiconvexity} in the second set of variables.
The notion of {\it curl Young quasiconvexity} has been introduced in \cite{APESAIM}, as a sufficient condition for the sequential weak *  lower semicontinuity of functionals of the type $\esssup_{x \in \Omega} f(x, \nabla u(x))$ in $W^{1,\infty}(\Omega;\mathbb R^d)$, (see also \cite{CDP} for a similar notion suited for $L^p$ approximation of supremal functionals, and \cite{RZ2} for the setting adopted in this paper).

\medskip
 
Let $Q$ be the unit cube $]0,1[^N$, and let $f:\mathbb{R}^{d \times N}\longrightarrow \mathbb{R}$ be a lower semicontinuous function, bounded from below. $f$ is {\it curl Young quasiconvex}, if 
$$f\left(\int_{\mathbb R^{d \times N}}\xi\, d\nu_x(\xi) \right)\leq \operatorname*{ess\,sup}_{y \in Q}  \left(\mathop{\nu_y - \rm ess\,sup}_{\xi \in \mathbb R^{d\times N}}f(\xi)\right),\ \mathcal L^N-\mathrm{a.e.}\ x\in Q,$$
whenever $\nu\equiv\{\nu_x\}_{x \in Q}$ is a $W^{1,\infty}$-gradient Young measure (see \cite{KP1991} for the introduction, \cite{Rindler} and \cite{Fonseca-Leoni-unpub} for a comprehensive description). For the readers' convenience we just say that 
	Young measures encode information on the oscillation behaviour of weakly converging sequences. 
	For a more detailed introduction to the topic, we refer to the broad literature, e.g.~\cite[Chapter~8]{Fonseca-Leoni-Lp}, \cite{Pedregal}, \cite[Section~4]{Rindler}.
	
\bigskip 

In general dimensions $l, m, n \in \mathbb N$, we denote by $\mathcal M(\mathbb R^l)$ the set of bounded Radon measures and by $\mathcal Pr(\mathbb R^l)$ its subsets of probability measures.

Let $U\subset \mathbb R^n$ be a Lebesgue measurable set with finite measure. 
By definition, a Young measure $\nu=\{\nu_x\}_{x\in U}$ is an element of the space $L^\infty_w(U;\mathcal M(\R^m))$ of essentially bounded, weakly$^\ast$ measurable maps defined in $U\to \mathcal M(\R^m)$, which is isometrically isomorphic to the dual of $L^1(U;C_0(\R^m))$, such that $\nu_x:=\nu(x) \in \mathcal Pr(\R^m)$ for $\mathcal L^n$-a.e.~$x \in U$.  
One calls $\nu$ 
homogeneous if there is a measure $\nu_0 \in  \mathcal Pr(\mathbb R^m)$ such that $\nu_x = \nu_0$ for $\mathcal L^n$- a.e.~$x \in U$.

A sequence $(z_j)_j$ of measurable functions $z_j: U\to \R^m$ is said to generate a Young measure $\nu \in L^\infty_w(U;\mathcal Pr(\R^m))$ 
if for every $h \in L^1(U)$ and $\varphi \in C_0(\R^m)$,
$$
\lim_{j\to\infty}\int_U h(x)\varphi(z_j(x))\,dx = \int_U h(x)\int_{\R^m}\varphi(\xi)d \nu_x(\xi)\,dx = \int_U h(x) \langle \nu_x, \varphi \rangle\,dx, 
$$  
or $\varphi(z_{j})\overset{\ast}{\rightharpoonup} \langle \nu_x,\varphi\rangle$ for all $\varphi \in C_0(\mathbb R^m)$; 
in formulas,
\begin{align*}
	z_j\stackrel{YM}{\longrightarrow} \nu \qquad \text{as $j\to \infty$,}
\end{align*} 
with $\langle \cdot, \cdot \rangle$ denoting the duality product between probability measures and continuous functions $C_0(\mathbb R^m)$.

To keep the brevity of this article we omit the fundamental theorem for Young measures, for which we refer to \cite[Theorems~8.2 and~8.6]{Fonseca-Leoni-Lp}, \cite[Theorem~4.1, Proposition~4.6]{Rindler}.

We also recall that if $(z_j)_j \subset L^p(U;\R^m)$ ($p \in (1,+\infty]$) generates a Young measure $\nu$
and converges weakly($^\ast$) in $L^p(U;\R^m)$  to a limit function $u$, then $[\nu_x] = \langle \nu_x, {\rm id}\rangle= \int_{\mathbb R^m}\xi d \nu_x(\xi) = u(x)$ for $\mathcal L^n$-a.e.~$x\in U$.  
In the sequel we will mainly restrict to gradient Young measure, namely with $U:=\Omega\subset \mathbb R^N$ a bounded open set, and $m=N\times d$, a $W^{1,\infty}$-gradient Young measure (see \cite{KP1991}) is a Young measure generated by a sequence of $(\nabla u_j)_j$ with $u_j \in W^{1,\infty}(\Omega;\mathbb R^d)$.

	For our purposes, we also recall that in \cite[Proposition 4.3 and Remark 4.4]{RZ2} {\it curl Young quasiconvexity} has been characterized as follows.

$f$  is {\it  curl Young quasiconvex} if and only if it verifies
\begin{align*}
	f\left(\int_{\mathbb R^{d \times N}}\xi\,d\nu(\xi)\right) \leq  \mathop{\nu - \rm ess\,sup}_{\xi \in \mathbb R^{d \times N}} f(\xi) 
\end{align*} 
whenever $\nu$ is a $W^{1,\infty}$-gradient Young measure.

\section{Lower semicontinuity}

The notion which will play a crucial role for us is 
the {\it separate curl Young-quasiconvexity}.

\begin{definition}\label{def}
	Let $W:\mathbb R^{d \times N}\times \mathbb R^{d \times N} \to \mathbb R$ be a lower semicontinuous function. $W$ is said to be {\it separately curl Young quasiconvex} if 
		\begin{align}
		W([\nu], [\mu])  & \leq  (\nu\otimes\mu)\text{-}\esssup_{(\xi, \zeta)\in \R^{d \times N}\times \R^{d \times N}} W(\xi, \zeta) \nonumber\\ &= \nu\text{-}\esssup_{\xi \in \R^{d \times N}} \bigl(\mu\text{-}\esssup_{\zeta \in \R^{d \times N}} W \xi, \zeta)\bigr) \label{WscYqcx}\\
		&= \mu\text{-}\esssup_{\zeta \in \R^{d \times N}} \bigl(\nu\text{-}\esssup_{\xi \in \R^{d \times N}} W(\xi, \zeta)\bigr) \nonumber
	\end{align}
 for every $\nu, \mu$ $ W^{1,\infty}$- gradient Young measures.
If $W:\Omega \times \Omega \times \R^{d \times N}\times \R^{d \times N}\to \R$ is a normal integrand bounded from below, then it
 is said to be   {\it separately curl Young quasiconvex} if $W(x,y, \cdot ,\cdot)$  is {\it separately curl Young quasiconvex } for $\mathcal L^N \otimes \mathcal L^N $-a.e. $(x,y)\in \Omega \times \Omega$. 
\end{definition}

A key tool for the proof of our result is the following lemma, first stated in \cite{B1} in the continuous and homogeneous case, and, then  proved in its current version in \cite{RZ2}.

\begin{lemma}\label{barronmon}
	Let $U\subset \mathbb R^n$ be an open set with finite measure and let $f:U\times \R^m\to \R$  be a normal integrand bounded from below. 		Further, let $(u_k)$ be a uniformly bounded sequence of functions in $L^\infty(U;\mathbb R^m)$ generating a Young measure $\nu=\{\nu_x\}_{x\in U}$. Then,
	\begin{equation*}
		\liminf_{k \to \infty}\operatorname*{ess\,sup}_{x\in U} f(x,u_k(x))  \geq \operatorname*{ess\,sup}_{x\in U}\bar f(x),
	\end{equation*}
	where $\bar f(x) := \operatorname*{\nu_{x}-ess\,sup}_{\xi\in \R^m} f(x,\xi)$ for $x\in U$.
\end{lemma}

%
%
%


With the aim of analyzing nonlocal problems, in \cite{KZ} to any function $u\in L^1(\Omega;\R^m)$ it has been associated the vector field 
\begin{align}\label{vu}
	w_u(x,y):=(u(x), u(y)) \quad \text{for $(x,y)\in \Omega\times \Omega$.}
\end{align}

In the sequel we will consider nonlocal fields $w_{\nabla u}(x,y)= (\nabla u(x), \nabla u(y))$ for $(x,y)\in \Omega\times \Omega$.

\medskip
	
The following lemma, which was established by Pedregal in~\cite[Proposition 2.3]{Pedregal}, gives a characterization of Young measures generated by sequences as in \eqref{vu}. 
\begin{lemma}
	\label{asprop2.3Pedregal}
	Let $(u_j)_j \subset L^p(\Omega;\R^m)$ with $1 \leq p\leq \infty$ generate a Young measure $\nu=\{\nu_x\}_{x\in \Omega}$, and let $\Lambda = \{\Lambda_{(x,y)}\}_{(x, y)\in \Omega\times \Omega}$ be a family of probability measures on $\mathbb R^m \times \mathbb R^m$. 
	
	Then $\Lambda$ is the Young measure generated by the sequence $(w_{u_j})_j\subset L^p(\Omega\times \Omega;\R^m\times \R^m)$ defined according to~\eqref{vu} 
	if and only if
	\begin{equation*} 
		\Lambda_{(x,y)} = \nu_x \otimes \nu_y \qquad \text{for $\mathcal L^N\otimes \mathcal L^N$-
			a.e.~$(x, y)\in \Omega\times \Omega$\color{black}}\end{equation*}
	and
	\begin{align*}
		\begin{cases}
			\displaystyle \int_{\Omega }\int_{\mathbb R^m} |\xi|^p \,d \nu_x(\xi)\,dx< \infty, & \hbox { if }p<\infty,\\[0.2cm]
			{\rm supp}\,\nu_x \subset  K  \hbox{ for $ \mathcal L^N$-a.e.~}x \in \Omega \hbox{ with a fixed compact set $K\subset \R^m$}, & \hbox{ if } p = \infty.
		\end{cases}
	\end{align*}
\end{lemma}

\begin{remark}
	\label{levelcx}
	The class of {\it separately  curl Young quasiconvexity} is not empty since any separately level convex function is {\it separately  curl Young
quasiconvex}, indeed in \cite[Lemma 3.5 (iv)]{KZ} it has been proven that any Borel function $W$, whose sublevel sets are separately convex (i.e. $W$ is separately level convex) satisfies
\eqref{WscYqcx} for every $\nu, \mu\in \mathcal Pr(\mathbb R^{d \times N})$. On the other hand the notions are not equivalent as we can see considering the function $W:\mathbb R^{2 \times 2}\times \mathbb R^{2 \times 2} \to \mathbb [0,+\infty]$, defined as 
$$W(\xi, \eta)= (\sup\{h(|\xi|), k(\xi)\} ) (\sup\{h(|\eta|), k(\eta)\} ),$$ with $h$ and $k$ as in \cite[Example 6.7]{APESAIM},  namely
$k(\Sigma):= \arctan (\rm det \Sigma)$ and

\noindent $h(t)= \left\{
\begin{array}{ll} 0 &\hbox{ if } t \leq 1,
\\
t-1 &\hbox{ if } 1 \leq  t \leq 2,\\
1 &\hbox{ if  }t \geq 2\end{array}
\right. $
Indeed for any fixed $\eta$ or $\xi$ the function $W(\cdot,\eta)$ or $W(\xi, \cdot)$ turns out to be {\it curl Young quasiconvex} but not  generally level convex.
\end{remark}

We are in position to establish our main result
\begin{theorem}
	\label{LsccurlYoung}
	Let $W:\Omega \times \Omega \times \mathbb R^{d \times N} \times \mathbb R^{d\times N} \to \mathbb R$ be a normal integrand, bounded from below and such that  $W(x,y,\cdot,\cdot)$ is separately {\it curl Young quasiconvex} for $\mathcal L^N \otimes \mathcal L^N$- a.e. $(x,y) \in \Omega\times \Omega$. Let $F:W^{1,\infty}(\Omega;\mathbb R^d)\to\mathbb R$ be the functional defined by 
	\begin{equation}\label{Fdef}
	F(u)= \esssup_{(x,y) \in \Omega\times \Omega }W(x,y,\nabla u(x),\nabla u(y)).
	\end{equation}
	Then the functional $F$ is sequentially weakly * lower semicontinuous in $W^{1,\infty}(\Omega;\mathbb R^d)$.
\end{theorem}

\begin{remark}\label{Aqcxlsc}
	We observe that this result extends to the non-homogeneous and differential setting \cite[Proposition 3.6]{KZ}. 
	
	The same proof could be used to show that separate level convexity of $ W(x,y,\cdot, \cdot)$ for $\mathcal L^N \otimes \mathcal L^N$-a.e. $(x,y) \in \Omega \times \Omega$ is sufficient to guarantee the sequential weak* lower semicontinuity  in $L^\infty(\Omega;\mathbb R^d)$ of $\esssup_{(x,y)\in \Omega \times \Omega }W(x,y, u(x), u(y))$.
	
	Nevertheless, as proven in the latter setting, under homogeneity assumptions, we conjecture that {\it separate curl Young quasiconvexity } is not `really' necessary for the sequential lower semicontinuity of the functional in \eqref{Fdef}, since from one hand some symmetry of $W$ should be taken into account, (cf. the notions of Cartesian separate level convexity in \cite{KZ, KRZ}), but also it is worth to observe that even in the local setting it is currently an open question the necessity of {\it curl Young quasiconvexity} for the sequential weak* lower semicontinuity of $\esssup_{x\in \Omega} f(\nabla u(x))$, namely it is not known, in general, if {\it curl Young quasiconvexity} is equivalent to  the {\it Strong Morrey quasiconvexity} introduced in \cite{BJW}, except some particular case as those considered in \cite{APESAIM} and \cite{PZ}.
	
	Finally, we also point out that, under suitable continuity conditions on the second set of variables for $W$, our arguments could be successfully employed to prove the lower semicontinuity of nonlocal supremal functionals under more general differential constraints than {\it curl}.
\end{remark}
\begin{proof}
	The result follows from Lemmas \ref{barronmon} \ref{asprop2.3Pedregal}, and Definition \ref{def}.
	Without loss of generality we can assume that $W$ is non negative.
	

		Let $(w_{\nabla u_j})_j\subset W^{1,\infty}(\Omega\times \Omega;\R^{d \times N}\times \R^{d \times N})$ be the sequence of nonlocal vector fields associated with $(\nabla u_j)_j$, cf.~\eqref{vu}, and $\Lambda=\{\Lambda\}_{(x,y)\in \Omega\times \Omega}=\nu_x\otimes \nu_y$ for $x,y\in \Omega$ the generated  $W^{1,\infty}$- gradient Young measure according to  Lemma~\ref{asprop2.3Pedregal}. 
	Lemma \ref{barronmon} implies that
	\begin{align}
		\nonumber
		\liminf_{j\to \infty} F(u_j) &= \liminf_{j\to \infty} \esssup_{(x,y)\in \Omega\times \Omega}W(x,y, \nabla u_j(x), \nabla u_j(y))\\ 
		&\geq \esssup_{(x,y)\in \Omega\times \Omega} \overline W(x,y), \label{one}
	\end{align}
	where $\overline W(x,y) := \Lambda_{(x,y)}\text{-}\esssup_{(\xi,\zeta)\in \mathbb R^{d \times N} \times \mathbb R^{d \times N}} W(\xi, \zeta)$. By Lemma \ref{asprop2.3Pedregal}, 
	\begin{align*}
		\overline W(x,y) &=  \nu_x\otimes \nu_y\text{-}\esssup_{(\xi,\zeta)\in \mathbb R^{d \times N} \times \mathbb R^{d \times N}} W(\xi, \zeta)\\
		&=\nu_x\text{-}\esssup_{\xi\in \mathbb R^{d \times N}}\bigl(\nu_y\text{-} \esssup_{\zeta \in \mathbb R^{d \times N}} W(\xi, \zeta) \bigr)
	\end{align*}
	for $\mathcal L^N \otimes \mathcal L^N$-a.e.~$(x, y)\in \Omega\times \Omega$, and since $W$ is separately {\it curl Young quasiconvex}, 
	
	it results that
	\begin{equation}
		\label{two}
		\overline W(x,y)
		\geq W(x,y, [\nu_x],[\nu_y]) = W(x,y, \nabla u(x), \nabla u(y))
	\end{equation}
	for $\mathcal L^N \otimes \mathcal L^N$-a.e. $(x,y) \in \Omega \times \Omega$.
	The proof follows from \eqref{one} and \eqref{two}.
	
\end{proof}

\section{Conclusions} In this paper we provide a sufficient condition for the lower semicontinuity of nonlocal supremal functionals depending on the gradients of suitable Lipschitz fields.  
We conjecture that this notion is also suitable to provide an $L^p$ approximation result in the spirit of what is proven for $L^\infty$ fields in \cite{KRZ}. This latter study and the search for necessary conditions will be the subject of future research.
We conclude observing that analogous results in the case of nonlocal integral functionals, depending on the gradient of scalar fields, can be found in \cite{Bevan Pedregal}.

\section*{Acknowledgements}
GG and EZ gratefully acknowledge support from INdAM GNAMPA.

The authors do not have any conflict of interest with third parties and they both contribute in the same way in the redaction of the manuscript.

Data availability statement: non applicable.

	\end{document}